\documentclass[10pt,a4paper]{article}
\usepackage[utf8]{inputenc}
\usepackage[english]{babel}
\usepackage{amsmath}
\usepackage{amsfonts}
\usepackage{amssymb}
\usepackage{hyperref}
\date{}
\author{Lukas Prader}
\title{A local-global principle \\ for surjective polynomial maps}
\begin{document}

\maketitle
\footnotetext{The author has been supported in part by project 
P-25652 of the Austrian Science Fund (FWF) and by a Qualifying Hausdorff Scholarship of the Bonn International Graduate School of Mathematics (BIGS-M).  
}

\begin{abstract} Let $R$ be an affine domain of characteristic zero with finite quotients. We prove that a polynomial map over $R$ is surjective if and only if it is surjective over $\widehat{R_{\mathfrak{m}}}$, the completion of $R$ with respect to $\mathfrak{m}$, for every maximal ideal $\mathfrak{m} \subseteq R$. In fact, the completions $\widehat{R_{\mathfrak{m}}}$ may be replaced by 
arbitrary subrings containing $R$. We use this result to yield a characterization of surjective polynomial maps, and remark that there does not exist a similar principle for injective polynomial maps.
\end{abstract}

All rings in this paper are supposed to be commutative with unity. Let now $R$ be a ring and $n$ a positive integer, write $x := (x_1,...,x_n)$.
Then any $n$-tuple $f = (f_1,...,f_n) \in R[x]^n$ of polynomials over $R$ in $n$ variables defines a map
$$f \colon R^n \to R^n, \text{ whereby } a = (a_1,...,a_n) \in R^n \mapsto \big( f_1(a),...,f_n(a) \big) \in R^n,$$
a so-called \textit{polynomial map} over $R$. As usually, we call $f$ surjective if $f \big( R^n \big) = R^n$, and injective if $f(a) \neq f(b)$ holds for all $a,b \in R^n, a \neq b$.
Further, $f$ is called invertible if $R[f] = R[x]$. 
This is the case if and only if there exists $g \in R[x]^n$ so that $f(g(x)) = g(f(x)) = x$, and consequently $\det (Jf)$ is a unit in $R[x]$, i.e. $\det (Jf) \in R[x]^\times$, where $Jf$ denotes the Jacobian matrix of $f$. \\
Let $S$ be an $R$-algebra, then it is important to observe that any $f  \in R[x]^n$ also gives rise to a polynomial map $f \colon S^n \to S^n$. For example, one may choose for $S$ the localization $R_{\mathfrak{m}}$ or completion $\widehat{R_{\mathfrak{m}}}$ of $R$ with respect to a maximal ideal $\mathfrak{m} \subseteq R$. This allows us to study polynomial maps locally, which leads to the following local-global principle. \\
\\
\textbf{Theorem 2.2:} \textit{Let $R$ be an affine domain of characteristic zero with finite quotients, i.e., $R/\mathfrak{m}$ is finite for every maximal ideal $\mathfrak{m} \subseteq R$. \\
Further, let $f \in R[x]^n$. Then
$$f \colon R^n \to R^n \text{ is surjective}$$
if and only if
$$f \colon \widehat{R_{\mathfrak{m}}}^n \to \widehat{R_{\mathfrak{m}}}^n \text{ is surjective for every maximal ideal } \mathfrak{m} \subseteq R.$$} \\
\\
In fact, it turns out that $\widehat{R_{\mathfrak{m}}}$ may even be replaced by an arbitrary subring of $\widehat{R_{\mathfrak{m}}}$ containing $R$, e.g., by the localization $R_{\mathfrak{m}}$. This will be discussed in Corollary 2.3. 
The notion of affine domains is adopted from \cite{Modelth} and will be briefly discussed in the first section of this paper. For instance, any domain of characteristic zero that is finitely generated as a $\mathbb{Z}$-algebra satisfies the requirements of Theorem 2.2. \\
We shall note that there exist results of a similar spirit in the literature. \\
For instance, it is a consequence of  \cite[Ex. 1.0.5,6]{BS} that $a \in \mathbb{Z}^n$ lies in the image of some $\mathbb{Z}$-linear map $f \colon \mathbb{Z}^n \to \mathbb{Z}^n$ if and only if $a$ lies in the image of the induced map $f \colon {\mathbb{Z}_p}^n \to {\mathbb{Z}_p}^n$ over the ring $\mathbb{Z}_p$ of $p$-adic integers for every prime $p$. In particular, this implies that Theorem 2.2 holds for maps of this kind. However, the stronger statement involving $a \in \mathbb{Z}^n$ fails for arbitrary polynomial maps, with $a= 0$ and $f = (p_1 x+1)(p_2x+1) \in \mathbb{Z}[x]$ for distinct primes $p_1,p_2$ as an obvious counterexample, see also Proposition 5.1.  \\
More generally, we shall mention Proposition 3.9 from \cite{AM}, asserting that a homomorphism $f \colon M \to N$ of $R$-modules is surjective (resp. injective) if and only if the induced homomorphism $f \colon M_{\mathfrak{m}} \to N_{\mathfrak{m}}$ is surjective (resp. injective) for every maximal ideal $\mathfrak{m} \subseteq R$. But indeed, we will see in the final section of this paper that both injectivity and non-injectivity of polynomial maps are no local properties. \\
It is a consequence of Theorem 2.2 that surjective polynomial maps are even invertible, and its main application is the following characterization of surjective resp. invertible polynomial maps, which (in the special case $R = \mathbb{Z}$) also appears in \cite[Thm. 10.3.13]{Arno1}. In fact, we are even able to drop the assumption on finite quotients here. \\
\\
\textbf{Proposition 3.2:} \textit{Let $R$ be an affine domain of characteristic zero but not a field, and let $f \in R[x]^n$. Then $f \colon R^n \to R^n$ is surjective
if and only if $\det(Jf) \in R^\times$ and the induced map $f \colon (R/\mathfrak{m})^n \to (R/\mathfrak{m})^n$ is bijective for every maximal ideal $\mathfrak{m} \subseteq R$.} \\
\\
The fourth section of this paper is devoted to the famous Jacobian problem. Indeed, we discuss its relation to Proposition 3.2 and provide an equivalent conjecture.

\section{Preparations}

We shall start with the following very powerful result, a proof can be found in \cite[Thm. 1.1.2]{Arno1}. \\
\\
\textbf{Result 1.1: (Formal inverse function theorem)} \textit{Let $R$ be a ring and let $f \in R[[x]]^n$ be such that $f(0) = 0$ and $\det((Jf)(0)) \in R^\times$. Then there exists a unique $g \in R[[x]]^n$ so that $g(0) = 0$ and $g(f(x)) = f(g(x)) = x$.}  \\
\\
It is clear that the invertibility of a polynomial map $f \colon R^n \to R^n$ implies the invertibility of $f \colon S^n \to S^n$ for any $R$-algebra $S$. Under further assumptions, Result 1.1 enables us to establish also the converse. \\
\\
\textbf{Corollary 1.2:} \textit{Let $R \subseteq S$ be an extension of rings and $f \in R[x]^n$. \\
(i) If $f$ is invertible over $S$ and $\det((Jf)(0)) \in R^\times$, then $f$ is invertible over $R$. \\
(ii) If $R[x]^\times = S[x]^\times \cap R[x]$, then the invertibility of $f \colon R^n \to R^n$ is equivalent to the invertibility of $f \colon S^n \to S^n$.} \\
\\
Note that if $f \colon S^n \to S^n$ is invertible and $R[x]^\times = S[x]^\times \cap R[x]$, then $\det (J f) \in S[x]^\times$ by the invertibility and $\det (J f) \in R[x]$ by assumption, so $\det (J f) \in R[x]^\times$ and $\det ((J f)(0)) \in R^\times$. Thus it suffices to deal with (i), a proof can be found in \cite[Lemma 1.1.8]{Arno1}. \\
\\
Let us briefly remark that the condition $R[x]^\times = S[x]^\times \cap R[x]$ is satisfied, for instance, if $L \mid K$ is an extension of number fields and $R = \mathrm{O}_K, S = \mathrm{O}_L$ are the rings of integers in $K,L$. This is due to the characterization of units in $\mathrm{O}_L$ and $\mathrm{O}_K$ by the respective field norms. \\
\\
Given a univariate polynomial over a domain $R$ without repeated zeros, we shall now discuss whether this also holds for the induced polynomial over some quotient $R/\mathfrak{m}$ by a maximal ideal. 
A proof of the forthcoming lemma can be found in \cite[Prop. 3.3.3]{Arno1}, it relies on basic facts about the resultant $\mathrm{Res}(f,g)$ of polynomials $f,g$.  \\
\\
\textbf{Lemma 1.3:} \textit{Let $R$ be a domain with field of fractions $K$, $\mathfrak{m} \subseteq R$ a maximal ideal and $f \in R[x]$ a univariate polynomial. \\
(i) $\mathrm{Res}(f,f') \neq 0$ if and only if $f$ does not have repeated zeros over $\overline{K}$, the algebraic closure of $K$. \\
(ii) If the resultant $\mathrm{Res}(f,f') \not\in \mathfrak{m}$, then the image of $f$ in $(R/\mathfrak{m})[x]$ does not have repeated zeros.} \\
\\
We shall continue with a couple of definitions, followed by a brief discussion. \\
A domain $R$ is called \textit{affine} if it is either finitely generated as a ring or finitely generated as an algebra over a subfield. In this paper, we mainly focus on affine domains of characteristic zero. Further, we say that an affine domain $R$ has \textit{finite quotients} if the quotient $R/\mathfrak{m}$ is finite for every maximal ideal $\mathfrak{m} \subseteq R$. For example, any domain $R$ that is finitely generated as a $\mathbb{Z}$-algebra is an affine domain with finite quotients, and we have $\mathrm{char}(R) = 0$ if and only if $\mathbb{Z} \subseteq R$. \\
Moreover, we say that a domain $R$ has \textit{enough nonunits} if for every univariate polynomial $f \in R[x] \setminus R$ we find $r \in R$ so that $f(r) \not\in R^\times$. Of course, any infinite domain $R$ with finite group of units $R^\times$ (e.g. $R = \mathbb{Z}$) has enough nonunits. On the other hand, such rings need to be semi-primitive, as otherwise $f = mx+1 \in R[x]$ with $0 \neq m \in Jac(R)$ has the property that $f(R) \subseteq R^\times$. \\
The following result can be found in \cite[2.2,2.9]{Modelth}. \\
\\
\textbf{Result 1.4:} \textit{(i) If a domain $R$ has enough nonunits, then for every $f \in R[x] \setminus R$ and $a \in R\setminus \{0\}$ there exist a maximal ideal $\mathfrak{m} \subseteq R$ and $r \in R$ so that $f(r) \in \mathfrak{m}$ but $a \not\in \mathfrak{m}$. 
(ii) Affine domains which are not fields have enough nonunits.} \\
\\
The proof of (i) is quite easy: If $f(0) = 0$, everything is clear. If not, one writes $f(a \cdot f(0) \cdot x) = f(0) \cdot h(x)$, where $h(x) = a \cdot x \cdot g(x) +1$ for some $g \in R[x]$. Since $R$ has enough nonunits, there exists $s \in R$ so that $h(s) \not\in R^\times$. In particular, we find a maximal ideal $\mathfrak{m} \subseteq R$ containing $h(s)$ and hence $f(r)$, where $r = a \cdot f(0) \cdot s$. But certainly $a \not\in \mathfrak{m}$, as otherwise $1 = h(r) - a \cdot r \cdot g(r) \in \mathfrak{m}$, which is impossible. \\
The second part can be proved by applying ideas from Galois theory. \\
\\
Our next intention is to discuss relations between injectivity and surjectivity of polynomial maps. The most familiar result seems to be the following, which is at least partially due to Ax and Grothendieck, see for instance \cite[Thm. 3.1]{Se}, also for further references. \\
\\
\textbf{Result 1.5:} \textit{Let $\mathbb{K}$ be an algebraically closed field, then any injective polynomial map $f  \colon \mathbb{K}^n \to \mathbb{K}^n$ is surjective. 
If further $\mathrm{char}(\mathbb{K}) = 0$, then $f$ is even invertible.} \\
\\
A detailed proof can be found in \cite[Thm.s 4.1.1, 4.2.1]{Arno1}. Concerning the first part, the idea is that the injectivity and the lack of surjectivity of $f$ could be expressed by polynomial identities, which is due to the Nullstellensatz. So one descends from $\mathbb{K}$ to a finite field in such a way that these identities remain valid and faces a contradiction, as surjectivity and injectivity are  equivalent over finite fields. The second part relies heavily on an algebro-geometric result about the fibres of morphisms, see \cite[Prop. B.2.1]{Arno1} or \cite[sect. 8]{Mumf}, which fails in prime characteristic. \\
Surprisingly, a reverse statement holds for affine domains. It is proved in \cite[Thm. A]{Modelth} and relies on the fact that affine domains have enough nonunits and are \textit{Hilbert dense} in their fields of fractions. The latter property involves techniques from mathematical logics. \\
\\
\textbf{Result 1.6:} \textit{Let $R$ be an affine domain but not a field, then any surjective polynomial map $f \colon R^n \to R^n$ is invertible, hence bijective.} \\
\\
It is clear that this argumentation via invertibility can only work in rings with the property that surjectivity and invertibility are equivalent for polynomial maps. For example, the ring $\mathbb{Z}_p$ does not have this property. To see this, note that the polynomial map induced by $f = px^2+x \in \mathbb{Z}_p[x]$ certainly is not invertible as $f' = 2px+1 \not\in \mathbb{Z}_p[x]^\times = \mathbb{Z}_p^\times$, but surjective by the prospective Proposition 2.1. \\
However, as the following proposition demonstrates, we can extend Result 1.6 to a different class of rings by avoiding the concept of invertibility. \\
\\ 
\textbf{Proposition 1.7:} \textit{Let $R$ be a ring containing a finitely generated ideal $\mathfrak{I}$ so that $R/\mathfrak{I}$ is finite and
$$\mathfrak{I}^\infty := \bigcap_{k=0}^\infty \frak{I}^k = \{0\}.$$
Then any surjective polynomial map $f \colon R^n \to R^n$ is bijective.} \\ 
\\
\textit{Proof:} The surjectivity of $f \colon R^n \to R^n$ clearly implies the surjectivity of the induced maps $f \colon (R/ \mathfrak{I}^k)^n \to (R/ \mathfrak{I}^k)^n$ for every $k \geq 1$. As $\mathfrak{I}$ is finitely generated and $R/\mathfrak{I}$ is finite, the quotients $R/ \mathfrak{I}^k$ are finite as well, and hence the induced maps are even bijective. Suppose now that $f(a) = f(b)$ for some $a,b \in R^n$, then\footnote{Given $x=(x_1,...,x_n), y=(y_1,...,y_n) \in R^n$, we say that $x \equiv y \text{ (mod } \mathfrak{m}^k \text{)}$ if $x_i \equiv y_i \text{ (mod } \mathfrak{m}^k \text{)}$ for every $1 \leq i \leq n$.} certainly $f(a) \equiv f(b) \text{ (mod } \mathfrak{I}^k \text{)}$ for every $k \geq 1$. The bijectivity of the induced maps now yields $a \equiv b \text{ (mod } \mathfrak{I}^k \text{)}$ for every $k \geq 1$, hence $a-b \in (\mathfrak{I}^\infty)^n = \{0\}^n$, implying $a=b$. This proves that $f$ is injective. $\hfill \Box$ \\
\\
\textbf{Remarks:} (i) The above proof basically shows that the induced map \\
$f \colon (R/ \mathfrak{I}^\infty)^n \to (R/ \mathfrak{I}^\infty)^n$ is injective. But we have $R/ \mathfrak{I}^\infty = R$ by assumption ($\mathfrak{I}^\infty = \{0\}$), implying that $f \colon R^n \to R^n$ is injective. \\
(ii) Let us briefly discuss to which kinds of rings this proposition can be applied. By Krull's intersection theorem, we know that $\mathfrak{I}^\infty = \{0\}$ for any proper ideal $\mathfrak{I}$ in a Noetherian domain or Noetherian local ring. Furthermore, it is a well-known fact that the quotient of a finitely generated algebra over $\mathbb{Z}$ or over a finite field with a maximal ideal is a finite field. Since the latter are also Noetherian by Hilbert's basis theorem, Proposition 1.7 can be applied to any finitely generated algebra over $\mathbb{Z}$ or over a finite field. In fact, it is also accessible to their completions with respect to maximal ideals: By \cite[Prop. 10.16, Thm. 10.26]{AM} these are Noetherian local rings, and the residue fields are isomorphic (hence of the same cardinality) by \cite[Prop. 10.15]{AM}. \\
(iii) To see that none of the conditions in Proposition 1.7 can be dropped in general, consider the ring $R = \mathbb{F}_2 \times \mathbb{C}$, where $\mathbb{F}_2$ denotes the field containing precisely two elements. The ideals of $R$ do not fulfil the requirements of the proposition, and in fact, the polynomial map defined by  
$$f = (0_{\mathbb{F}_2},1_{\mathbb{C}}) \cdot x^2 + (1_{\mathbb{F}_2},0_{\mathbb{C}}) \cdot x \in R[x]$$
is surjective by the fundamental theorem of algebra, but not injective as \\ $f(0_{\mathbb{F}_2},1_{\mathbb{C}}) = f(0_{\mathbb{F}_2},-1_{\mathbb{C}})$. \\
\\
Before we are able to prove the key lemma of this paper, we need to recall one more result from commutative algebra. It is a consequence of \cite[Thm. 73]{Mat}. \\
\\
\textbf{Result 1.8:} \textit{Let $R$ be a finitely generated domain over $\mathbb{Z}$ or a field. Then 
$$\mathrm{Sing}(\mathrm{Spec}(R)) = \{ \mathfrak{
p} \subseteq R \text{ prime} \mid R_\mathfrak{p} \text{ is not regular}\} \subsetneq \mathrm{Spec}(R)$$
is Zariski-closed, i.e., of the form $\mathrm{V}(\mathfrak{I})$ for some ideal $(0) \subsetneq \mathfrak{I} \subseteq R$.} \\
\\
Note that $\mathrm{Sing}(\mathrm{Spec}(R))$ is a proper subset of $\mathrm{Spec}(R)$ since $R$ is a domain. Indeed, this means that $(0) \in \mathrm{Spec}(R)$, and $R_{(0)}$, the field of fractions of $R$, is certainly regular. \\ 
\\
Now we are ready to establish the main ingredient for the proof of Theorem 2.2. Basically, it is a generalization of \cite[Thm. 10.3.1]{Arno1}. \\
\\
\textbf{Lemma 1.9:} \textit{Let $R$ be an affine domain but not a field, $\mathrm{char}(R)=0$ and $K$ its field of fractions. Then for any $a_1,...,a_m \in \overline{K}$ we find infinitely many maximal ideals $\mathfrak{m} \subseteq R$ so that there is an injective $R$-homomorphism of rings}
$$\phi \colon R[a_1,...,a_n] \to \widehat{R_{\mathfrak{m}}}.$$ \\
\textit{Proof:} As $\mathrm{char}(R) = \mathrm{char}(K) = 0$, we find $a \in \overline{K}$ so that $K[a_1,...,a_m] = K[a]$ by the primitive element theorem. Thus there exist $f_j \in K[x]$ so that $f_j(a) = a_j$ for every $1 \leq j \leq m$. Choose $s \in R\setminus \{0\}$ so that $sf_j \in R[x]$ for each $j$, then we have inclusion maps $R[a_1,...,a_m] \to R[s^{-1}][a] \to K[a]$. \\
Let $f$ be the minimal polynomial of $a$ over $K$. Multiplying $f$ by a non-zero element from $R$, we may assume that $f \in R[x]$. Since $\mathrm{char}(K) = 0$, we know that $f$ does not have repeated zeros in $\overline{K}$, thus $d := \mathrm{Res}(f,f') \neq 0$ by Lemma 1.3(i). Further, we know that $\mathrm{Sing}(\mathrm{Spec}(R)) = V(\mathfrak{I})$ for some ideal $(0) \neq \mathfrak{I} \subseteq R$ by Result 1.8. Let us choose $i \in \mathfrak{I}\setminus \{0\}$ arbitrarily. \\
By Result 1.4 we now find a maximal ideal $s \cdot d \cdot i \not\in \mathfrak{m} \subseteq R$ so that $f$ has a root modulo $\mathfrak{m}$. As $d \not\in \mathfrak{m}$, this is indeed a simple root by Lemma 1.3(ii), proving\footnote{One can proceed as in the proof of Proposition 2.1, $(ii) \Longrightarrow (i)$.} the existence of $b \in \widehat{R_{\mathfrak{m}}}$ so that $f(b)=0$ by Hensel's lifting lemma. Moreover, $i \not\in \mathfrak{m}$ implies that $\mathfrak{m} \not\in \mathrm{Sing}(\mathrm{Spec}(R))$, hence $R_{\mathfrak{m}}$ and $\widehat{R_{\mathfrak{m}}}$ are regular local rings, thus integral domains. This yields a $K$-homomorphism $K[a] \to \mathrm{Quot}(\widehat{R_{\mathfrak{m}}})$ with the property that $a \mapsto b$. Now define $\phi$ as the composition of injective $R$-homomorphisms
$$\phi \colon R[a_1,...,a_m] \to R[s^{-1}][a] \to K[a] \to \mathrm{Quot}(\widehat{R_{\mathfrak{m}}}).$$
As $s \not\in \mathfrak{m}$, the elements of $R[s^{-1}]$ map into $\widehat{R_{\mathfrak{m}}}$ under $R[s^{-1}][a] \to K[a] \to \mathrm{Quot}(\widehat{R_{\mathfrak{m}}})$. Since we know that $a$ maps to $b \in \widehat{R_{\mathfrak{m}}}$, this proves that the image of $R[s^{-1}][a]$ under $R[s^{-1}][a] \to K[a] \to \mathrm{Quot}(\widehat{R_{\mathfrak{m}}})$ is contained in $\widehat{R_{\mathfrak{m}}}$, and thus the same holds for $R[a_1,...,a_m]$, so $\phi$ has the desired property. To see that there are infinitely many possible choices for $\mathfrak{m}$, one can simply replace the element $s \cdot d \cdot i$ in the above argument by $\delta \cdot s \cdot d \cdot i$ for some $\delta \in \mathfrak{m} \setminus \{0\}$.  $\hfill \Box$ 

\section{The principle and its consequences}

As usual, we shall first deal with the local situation. \\
\\
\textbf{Proposition 2.1:} \textit{Let $R$ be a Noetherian local ring with maximal ideal $\mathfrak{m} \neq (0)$ so that $R/\mathfrak{m}$ is finite and $R$ is complete with respect to the $\mathfrak{m}$-adic topology. Further, let $f \in R[x]^n$. Then the following assertions are equivalent: \\
(i) 
$$f \colon R^n \to R^n \text{ is surjective}.$$
(ii)
$$f \colon (R/\mathfrak{m})^n \to (R/\mathfrak{m})^n \text{ is bijective and } \det((Jf)(a)) \in R^\times \text{ for every } a \in R^n.$$}
\\
\textit{Proof:} If $f \colon R^n \to R^n$ is surjective, then clearly $f \colon (R/\mathfrak{m}^k)^n \to (R/\mathfrak{m}^k)^n$ is surjective and hence bijective as $R/\mathfrak{m}^k$ is finite for every $k \geq 1$. Now suppose to the contrary that $\det((Jf)(a)) \in \mathfrak{m} = R \setminus R^\times$ for some $a \in R^n$. Thus $(Jf)(a)$ is not invertible over $R/ \mathfrak{m}$, so we find $x \in R^n$ so that $x \not\in \mathfrak{m}^n$ but $(Jf)(a) \cdot x \in \mathfrak{m}^n$. 
Further, note that for any $m \in \mathfrak{m}^n$ we have that
$$f(a+m) \equiv f(a) + (Jf)(a) \cdot m \text{ (mod } \mathfrak{m}^2 \text{)}.$$
In particular, this holds for $m = \rho \cdot x$, where $\rho \in \mathfrak{m} \setminus \mathfrak{m}^2$ is arbitrary (note that $\rho$ exists as $R$ is Noetherian local). Since
$$(Jf)(a) \cdot m = \rho \cdot (Jf)(a) \cdot x \equiv 0 \text{ (mod } \mathfrak{m}^2 \text{)},$$
hence $f(a+m) \equiv f(a) \text{ (mod } \mathfrak{m}^2 \text{)}$, but $m = \rho \cdot x \not\equiv 0 \text{ (mod } \mathfrak{m}^2 \text{)}$ by construction, we see that $f \colon (R/\mathfrak{m}^2)^n \to (R/\mathfrak{m}^2)^n$ fails to be injective, a contradiction. \\
The fact that (ii) implies (i) is basically due to Hensel's lifting lemma: Replacing $f$ by $f-b$ for arbitrary $b \in R^n$, it suffices to prove the existence of $a \in R^n$ so that $f(a) = 0$.
This can be done by inductively constructing a Cauchy sequence $(a_k)_{k \geq 1}$ with the property that $f(a_k) \equiv 0 \text{ (mod } \mathfrak{m}^k \text{)}$ for every $k \geq 1$. 
By the bijectivity of $f \colon (R/\mathfrak{m})^n \to (R/\mathfrak{m})^n$, we are certainly able to choose $a_1 \in R^n$ as desired. Given $k \geq 1$ and $a_k \in R^n$ so that $f(a_k) \equiv 0 \text{ (mod } \mathfrak{m}^k \text{)}$, we recall that $(Jf)(a_k)$ is invertible over $R$ by assumption and define  $a_{k+1}:=a_k - ((Jf)(a_k))^{-1} f(a_k)$. Then by Taylor's formula,
$$f(a_{k+1}) \equiv f(a_k) - (Jf)(a_k) \cdot ((Jf)(a_k))^{-1} f(a_k) = 0 \text{ (mod } \mathfrak{m}^{k+1} \text{)}.$$
Finally, note that $(a_k)_{k \geq 1}$ is Cauchy since $a_k \equiv a_l \text{ (mod } \mathfrak{m}^l \text{)}$ for every $k \geq l$, hence has a limit $a \in R^n$ which suffices $f(a) = 0$ by construction. $\hfill \Box$ \\
\\
Now we are finally prepared to prove our main result, the argumentation is inspired by \cite[Thm. 10.3.8]{Arno1}. \\
\\
\textbf{Theorem 2.2:} \textit{Let $R$ be an affine domain of characteristic zero with finite quotients, and let $f \in R[x]^n$. Then
$$f \colon R^n \to R^n \text{ is surjective}$$
if and only if
$$f \colon \widehat{R_{\mathfrak{m}}}^n \to \widehat{R_{\mathfrak{m}}}^n \text{ is surjective for every maximal ideal } \mathfrak{m} \subseteq R.$$} \\
\textit{Proof:} As $R / \mathfrak{m}$ is finite for every maximal ideal $\mathfrak{m}$, we know that $\widehat{R_{\mathfrak{m}}}$ is a compact topological space. Thus $f(\widehat{R_{\mathfrak{m}}}^n) \subseteq \widehat{R_{\mathfrak{m}}}^n$ is closed, and contains $R^n$ if $f \colon R^n \to R^n$ is surjective. As the latter is a dense subset of $\widehat{R_{\mathfrak{m}}}^n$, this proves that $f(\widehat{R_{\mathfrak{m}}}^n) = \widehat{R_{\mathfrak{m}}}^n$, thus $f \colon \widehat{R_{\mathfrak{m}}}^n \to \widehat{R_{\mathfrak{m}}}^n$ is surjective for every $\mathfrak{m} \subseteq R$ maximal. \\
Suppose now that $f \colon \widehat{R_{\mathfrak{m}}}^n \to \widehat{R_{\mathfrak{m}}}^n$ is surjective (and hence injective by Proposition 1.7) for every maximal ideal $\mathfrak{m} \subseteq R$. We claim that this implies the injectivity of $f \colon \overline{K}^n \to \overline{K}^n$, where $K$ denotes the field of fractions of $R$. To see this, assume that $f(a) = f(b)$ for some $a,b \in \overline{K}^n$. By Lemma 1.9, we find a maximal ideal $\mathfrak{m} \subseteq R$ and an injective $R$-homomorphism $\phi \colon R[a,b] \to \widehat{R_{\mathfrak{m}}}$. Denoting the extended map $R[a,b]^n \to \widehat{R_{\mathfrak{m}}}^n$ also by $\phi$, we have that  
$$f(\phi(a)) = \phi(f(a)) = \phi(f(b)) = f(\phi(b)),$$
so $\phi(a) = \phi(b)$ by the injectivity of $f \colon \widehat{R_{\mathfrak{m}}}^n \to \widehat{R_{\mathfrak{m}}}^n$ and thus $a = b$ by the injectivity of $\phi$. So $f \colon \overline{K}^n \to \overline{K}^n$ is injective and thus invertible over $\overline{K}$ by Result 1.5, implying that $\det(Jf) \in R[x] \cap \overline{K}^\times = R\setminus \{0\}$.  But Proposition 2.1 tells us that the surjectivity of $f \colon \widehat{R_{\mathfrak{m}}}^n \to \widehat{R_{\mathfrak{m}}}^n$ implies that $\det((Jf)(a)) \in R \cap (\widehat{R_{\mathfrak{m}}})^\times = R \setminus \mathfrak{m}$ for every $a \in R^n$. Because this holds for every maximal ideal $\mathfrak{m} \subseteq R$, we get that
$$\det(Jf) \in \bigcap_{\mathfrak{m} \text{ max.}} (R \setminus \mathfrak{m}) = R \setminus \big( \bigcup_{\mathfrak{m} \text{ max.}} \mathfrak{m} \big) = R^\times.$$
As a consequence, Corollary 1.2(i) yields that $f \colon R^n \to R^n$ is invertible, hence surjective. $\hfill \Box$ \\
\\
Note that the surjectivity of $f \colon \widehat{R_{\mathfrak{m}}}^n \to \widehat{R_{\mathfrak{m}}}^n$ for all but finitely many maximal ideals $\mathfrak{m}$ is sufficient for the invertibility of $f \colon \overline{K}^n \to \overline{K}^n$, but not for the invertibility of $f \colon R^n \to R^n$. As an easy example, consider $f = 2x \in \mathbb{Z}[x]$. Then $f \colon \mathbb{Z}_p \to \mathbb{Z}_p$ is surjective for every prime $p \neq 2$, but $f \colon \mathbb{Z}_2 \to \mathbb{Z}_2$ and $f \colon \mathbb{Z} \to \mathbb{Z}$ fail to be surjective as $2$ is not a unit in these rings. Nevertheless, $f \colon \overline{\mathbb{Q}} \to \overline{\mathbb{Q}}$ is invertible with inverse given by $x \mapsto \frac{1}{2}x$, which however does not define a polynomial map over $\mathbb{Z}$. \\
\\
Given the theorem, we can prove a more general result without any effort. \\
\\
\textbf{Corollary 2.3:} \textit{Let $R$ be an affine domain of characteristic zero with finite quotients, and let $f \in R[x]^n$. For every maximal ideal $\mathfrak{m} \subseteq R$ we choose a subring $S(\mathfrak{m})$ of $\widehat{R_{\mathfrak{m}}}$ containing $R$. Then
$$f \colon R^n \to R^n \text{ is surjective}$$
if and only if 
$$f \colon S(\mathfrak{m})^n \to S(\mathfrak{m})^n \text{ is surjective for every maximal ideal } \mathfrak{m} \subseteq R.$$}
\textit{Proof:} It is a consequence of Theorem 2.2 that any surjective polynomial map $f \colon R^n \to R^n$ must be invertible, implying the (invertibility and hence) surjectivity of $f \colon S(\mathfrak{m})^n \to S(\mathfrak{m})^n$ for every $\mathfrak{m}$. On the other hand, $S(\mathfrak{m})$ must be a dense subset of $\widehat{R_{\mathfrak{m}}}$ as it contains $R$. That means, for every $\mathfrak{m} \subseteq R$ maximal, the surjectivity of $f \colon S(\mathfrak{m})^n \to S(\mathfrak{m})^n$ implies the surjectivity of $f \colon \widehat{R_{\mathfrak{m}}}^n \to \widehat{R_{\mathfrak{m}}}^n$. But now Theorem 2.2 yields the claim. $\hfill \Box$ \\
\\
\textbf{Remark:} Let $R$ and $f$ be as above. The following assertions are either equivalent versions or special cases of Theorem 2.2 and Corollary 2.3. \\
(i) $f \colon R^n \to R^n$ is surjective if and only if $f \colon \widehat{R_{\mathfrak{p}}}^n \to \widehat{R_{\mathfrak{p}}}^n$ is surjective for every prime ideal $\mathfrak{p} \subseteq R$. \\
(ii) $f \colon R^n \to R^n$ is $\blacktriangle$ if and only if $f \colon \widehat{R_{\mathfrak{m}}}^n \to \widehat{R_{\mathfrak{m}}}^n$ is $\blacktriangledown$ for every maximal ideal $\mathfrak{m} \subseteq R$. Here, one may choose $\blacktriangle, \blacktriangledown \in \{ \text{surjective, bijective, invertible} \}$ arbitrarily. \\
(iii) $f \colon R^n \to R^n$ is surjective if and only if $f \colon {R_{\mathfrak{m}}}^n \to {R_{\mathfrak{m}}}^n$ is surjective for every maximal ideal $\mathfrak{m} \subseteq R$. \\
(iv) $f \colon R^n \to R^n$ is surjective if and only if $f \colon R[s^{-1}]^n \to R[s^{-1}]^n$ is surjective for every nonunit $s \in R \setminus \{0\}$. \\
(v) Let $p$ be a prime number. Then $f \colon \mathbb{Z}^n \to \mathbb{Z}^n$ is surjective if and only if $f \colon {\mathbb{Z}[p^{-1}]}^n \to {\mathbb{Z}[p^{-1}]}^n$ and $f \colon {\mathbb{Z}_{(p)}}^n \to {\mathbb{Z}_{(p)}}^n$ are surjective. \\
\\
One of the advantages of the above local-global principle is that it immediately yields a precise characterization of surjective and invertible polynomial maps. \\
\\
\textbf{Corollary 2.4:} \textit{Let $R$ and $f$ be as above, then the following assertions are equivalent: \\
(i) $f \colon R^n \to R^n$ is surjective. \\
(ii) $f \colon R^n \to R^n$ is invertible. \\
(iii) $\det(Jf) \in R^\times$ and the induced map $f \colon (R/\mathfrak{m})^n \to (R/\mathfrak{m})^n$ is bijective for every maximal ideal $\mathfrak{m} \subseteq R$.} \\
\\
Note that the implication (iii) $\Longrightarrow$ (i),(ii) is due to Proposition 2.1. Moreover, if $R = \mathrm{O}_K$ is the ring of integers of some number field $K$, then Corollary 1.2(ii) and its subsequent remark allow us to add the further assertions \\
\\
\textit{(iv) $f \colon \mathrm{O}_L^n \to \mathrm{O}_L^n$ is surjective (or equivalently, invertible) for some finite field extension $L \mid K$. \\
(v) $f \colon \mathrm{O}_L^n \to \mathrm{O}_L^n$ is surjective (or equivalently, invertible) for any finite field extension $L \mid K$.} \\
\\
Finally, let us recall that we used in the proof of the theorem that $f \colon \widehat{R_{\mathfrak{m}}}^n \to \widehat{R_{\mathfrak{m}}}^n$ is surjective if and only if $R^n \subseteq f( \widehat{R_{\mathfrak{m}}}^n )$. We shall apply this fact in order to rephrase Theorem 2.2 in terms of Diophantine equations. \\
\\
\textbf{Corollary 2.5:} \textit{Let $R$ and $f$ be as above. Write $f = (f_1,...,f_n)$ and consider the simultaneous equations 
$$f_1(x_1,...,x_n) = a_1$$
$$\vdots$$
$$f_n(x_1,...,x_n) = a_n$$
for all $(a_1,...,a_n) \in R^n$. \\
Then each system possesses a solution in $R^n$ if and only if each system possesses a solution in $\widehat{R_{\mathfrak{m}}}^n$ for every $\mathfrak{m} \subseteq R$ maximal.}

\section{Dropping the assumption of finite quotients}

Let us return once again to the proof of Proposition 2.1. While the assumption on $R/\mathfrak{m}$ to be finite was not necessary for (ii)$\Longrightarrow$(i), the proof of the other implication certainly collapses if we drop that assumption. In fact, as the following example illustrates, finite quotients are even crucial at this point. \\
\\
\textbf{Example 3.1:} Consider the polynomial $f = x^3-3x \in \mathbb{C}[t][x]$, we shall discuss some of its properties: \\
(i) The map $f \colon \mathbb{C}[t] \to \mathbb{C}[t]$ is neither injective nor surjective. To see this, note that $f(0) = f(\sqrt{3})$, and that $f(a) = t$ for some $a \in \mathbb{C}[t]$ would imply that $\deg(a)<1$, hence $a \in \mathbb{C}$, which is impossible as $f( \mathbb{C}) \subseteq \mathbb{C}$. \\
(ii) In particular, the induced map $f \colon \mathbb{C}[t]/(t) \to \mathbb{C}[t]/(t)$ is not injective. However, it is surjective since $\mathbb{C}[t]/(t) \cong \mathbb{C}$ is algebraically closed. \\
(iii) Nevertheless, the map $f \colon \widehat{\mathbb{C}[t]_{\mathfrak{m}}} \to \widehat{\mathbb{C}[t]_{\mathfrak{m}}}$ is surjective for every maximal ideal $\mathfrak{m} \subseteq \mathbb{C}[t]$. To give a proof, we first note that by the Nullstellensatz there exists $\lambda \in \mathbb{C}$ such that 
$\mathfrak{m} = (t-\lambda)\mathbb{C}[t]$. Without loss of generality, we may assume that $\lambda = 0$, hence that $\widehat{\mathbb{C}[t]_{\mathfrak{m}}} \cong \mathbb{C}[[t]].$ Now let $s = \sum_{k=0}^\infty c_k t^k \in \mathbb{C}[[t]]$ be arbitrary, we want to construct $p \in \mathbb{C}[[t]]$ so that $f(p)=s$. By Hensel's lifting lemma and the fact that $\mathbb{C}[[t]]/(t) \cong \mathbb{C}[t]/(t) \cong \mathbb{C}$, it suffices to find $a \in \mathbb{C}$ such that the following holds:
$$f(a) = a^3-3a =c_0,$$
$$f'(a) = 3a^2-3 \neq 0.$$
Since $\mathbb{C}$ is algebraically closed, the above equality is  solvable for every $c_0$. But the inequality suggests that $a \not\in \{-1,1\}$, which is certainly the case if $c_0 \not\in \{-2,2\}$. However, if $c_0 =  \pm 2$ we may choose $a= \pm 2$, which shows that our desired $a \in \mathbb{C}$ always exists. \\
\\
As a bottom line, this provides a counterexample to Proposition 1.7, Proposition 2.1 and Theorem 2.2 if the assumption of finite quotients is dropped. In particular, Example 3.1 reveals that in the setting of Theorem 2.2 (with possibly infinite quotients), the surjectivity of $f \colon \widehat{R_{\mathfrak{m}}}^n \to \widehat{R_{\mathfrak{m}}}^n$ for every $\mathfrak{m} \subseteq R$ maximal is not strong enough to force the surjectivity of $f \colon R^n \to R^n$. \\ 
However, the assumption on finite quotients in Theorem 2.2 can be dropped if we require the maps $f \colon \widehat{R_{\mathfrak{m}}}^n \to \widehat{R_{\mathfrak{m}}}^n$ to be invertible. This will appear as a by-product of the following generalization of Corollary 2.4. \\
\\
\textbf{Proposition 3.2:} \textit{Let $R$ be an affine domain of characteristic zero but not a field, and let $f \in R[x]^n$. Then $f \colon R^n \to R^n$ is surjective
if and only if $\det(Jf) \in R^\times$ and the induced map $f \colon (R/\mathfrak{m})^n \to (R/\mathfrak{m})^n$ is bijective for every maximal ideal $\mathfrak{m} \subseteq R$.} \\
\\
\textit{Proof:} If $f \colon R^n \to R^n$ is surjective, then it is also invertible by Result 1.6 and hence the claim follows. \\
On the other hand, if $\det(Jf) \in R^\times$ and the induced map $f \colon (R/\mathfrak{m})^n \to (R/\mathfrak{m})^n$ is bijective for some $\mathfrak{m} \subseteq R$ maximal, then this implies the bijectivity of $f \colon \widehat{R_{\mathfrak{m}}}^n \to \widehat{R_{\mathfrak{m}}}^n$. To see this, note that surjectivity follows from Proposition 2.1. And in fact, it is easy to convince oneself that the choice of the $(a_k)_{k \geq 1}$ is indeed unique, which is mainly due to the invertibility of $(Jf)(a)$ over $R$ for every $a \in R^n$. Thus we are in the situation of Theorem 2.2 (i.e., injectivity of the maps $f \colon \widehat{R_{\mathfrak{m}}}^n \to \widehat{R_{\mathfrak{m}}}^n$ and $\det(Jf) \in R^\times$), implying that $f \colon R^n \to R^n$ is surjective. $\hfill \Box$ 

\section{Direct limits and the Jacobian problem}

Let us now focus on the less general setting of subrings $R$ of $\mathbb{Q}$. Then certainly $\mathbb{Z} \subseteq R$ and we can write $R = S^{-1}\mathbb{Z}$ for some multiplicatively closed set $S \subseteq \mathbb{Z}$. In particular, $R$ is a principal ideal domain, and its prime elements are (up to multiplication by units) precisely those prime numbers $p \in \mathbb{Z}$ which \textit{do not meet} $S$, i.e., such that $p \mathbb{Z} \cap S = \emptyset$. If we further know that $R$ is finitely generated over $\mathbb{Z}$, then Corollary 2.4 asserts that a polynomial map $f \colon R^n \to R^n$ is invertible if and only if $\det(Jf) \in R^\times$ and $f \colon (R/pR)^n \to (R/pR)^n$ is bijective for every prime $p$ not meeting $S$. If $f$ satisfies the latter two properties, then we shall say that $f$ \textit{fulfils $\mathcal{P}(R)$}. Now it is interesting to discuss what can be said if $R$ is not finitely generated over $\mathbb{Z}$. \\
First note that it suffices to study invertible polynomial maps over $R$ induced by elements in $\mathbb{Z}[x]^n$, as $x \mapsto r \cdot x$ is certainly invertible for every $r \in R^\times$. \\
The key idea is to view $R$ as a direct limit $R = \varinjlim \mathbb{Z}[p^{-1}]$, where $p$ ranges through all prime numbers $p$ that meet $S$. In particular, this emphasizes that every polynomial map fulfilling $\mathcal{P}(\mathbb{Z}[p^{-1}])$ is also invertible over $R$. We now want to define $\varinjlim \mathcal{P}(\mathbb{Z}[p^{-1}])$ as the strongest condition containing $\mathcal{P}(\mathbb{Z}[p^{-1}])$ for every $p$. This can reasonably be done in the following way: Let $f \in \mathbb{Z}[x]^n$ fulfil $\varinjlim \mathcal{P}(\mathbb{Z}[p^{-1}])$ if $\det(Jf) \in R^\times$ and $f \colon (R/pR)^n \to (R/pR)^n$ is bijective for every prime $p$ not meeting $S$. This brings us to the forthcoming \\
\\
\textbf{Conjecture 4.1:} \textit{Let $(R_i)_{i \in I}$ be a family of subrings of $\mathbb{Q}$ finitely generated over $\mathbb{Z}$, and write $$R = \varinjlim_{i \in I} R_i.$$
Then $\varinjlim_{i \in I} \mathcal{P}(R_i)$ is a sufficient condition for the invertibility of polynomial maps over $R$.} \\
\\
Note that Conjecture 4.1 is certainly true if $R$ is a finitely generated $\mathbb{Z}$-algebra. If it was true in general, 
then we could easily  characterize the invertible polynomial maps of any subring $R$ of $\mathbb{Q}$: The polynomial map $f \colon R^n \to R^n$ would be invertible if and only if $\det(Jf) \in R^\times$ and $f \colon (R/pR)^n \to (R/pR)^n$ was bijective for every prime $p$ not meeting $S$. \\   
For example, if $R = \mathbb{Z}_{(p)}$ for some prime $p$, then the conjecture would imply that $f \colon {\mathbb{Z}_{(p)}}^n \to {\mathbb{Z}_{(p)}}^n$ is invertible if and only if $\det(Jf) \in {\mathbb{Z}_{(p)}}^\times$ and $f \colon (\mathbb{Z}_{(p)}/p\mathbb{Z}_{(p)})^n \to (\mathbb{Z}_{(p)}/p\mathbb{Z}_{(p)})^n$ is bijective.
But the most interesting case seems to be $R = \mathbb{Q}$. Here we  have $S = \mathbb{Z}\setminus \{0\}$, so every prime in $\mathbb{Z}$ meets $S$, and thus the conjecture simply states that $f \colon \mathbb{Q}^n \to \mathbb{Q}^n$ is invertible if $\det(Jf) \in \mathbb{Q}^\times$, which is precisely the Jacobian conjecture for $\mathbb{Q}$. In fact, even the following holds. \\
\\
\textbf{Proposition 4.2:} \textit{Conjecture 4.1 is equivalent to the Jacobian conjecture.} \\
\\
\textit{Proof:} If the Jacobian conjecture over $\mathbb{Q}$ was true, then it would also hold over any subring $R_i$ of $\mathbb{Q}$ by \cite[Thm. 1.1.18]{Arno1}. Thus the claim follows from the definition of $\varinjlim_{i \in I} \mathcal{P}(R_i)$. $\hfill \Box$ 

\section{Final comments on injective polynomial maps}

In this final section, we shall move our focus from surjective to injective polynomial maps. For the sake of simplicity, we  restrict to the case of polynomial maps over $\mathbb{Z}
$, and we will be able to illustrate all phenomena of our interest by means of 1-dimensional polynomial maps. \\
\\
Reflecting Theorem 2.2, we now might wonder whether an analogue result also holds for injective polynomial maps. However, the following proposition demonstrates that this is not the case. \\
\\
\textbf{Proposition 5.1:}  \textit{There are injective polynomial maps $f \colon \mathbb{Z} \to \mathbb{Z}$  with the property that $f \colon \mathbb{Z}_{(p)} \to \mathbb{Z}_{(p)}$ is not injective for every prime number $p$.} \\
\\
Since $\mathbb{Z}_{(p)}$ may be regarded as a subring of $\mathbb{Z}_p$, it is then immediately clear that the maps $f \colon \mathbb{Z}_p \to \mathbb{Z}_p$ are not injective as well. \\ 
\\
\textit{Proof:} Consider $f = (p_1x-1) \cdot (p_2x-1) \cdot x \in \mathbb{Z}[x]$, where $p_1 < p_2$ are prime numbers. Then $f \colon \mathbb{Z} \to \mathbb{Z}$ is injective as $f(0) < f(1)$ and $f'(x) > 0$ for all $x \in (- \infty, 0] \cup [\frac{1}{p_1},\infty)$. 
On the other hand, $f$ has two distinct zeros in $\mathbb{Z}_{(p)}$ for every $p$, namely $0$ and $\frac{1}{p_1}$ if $p \neq p_1$ resp. $0$ and $\frac{1}{p_2}$ if $p \neq p_2$. Therefore none of the maps $f \colon \mathbb{Z}_{(p)} \to \mathbb{Z}_{(p)}$ can be injective. $\hfill \Box$ \\
\\
Thus Proposition 5.1 proves that Theorem 2.2 does not hold for injective (instead of surjective) polynomial maps. In fact, we even see that the non-injectivity of polynomial maps is not a local property either, as we could otherwise find for any injective $f \colon \mathbb{Z} \to \mathbb{Z}$ a prime number $p$ so that 
$f \colon \mathbb{Z}_{(p)} \to \mathbb{Z}_{(p)}$ is injective. \\
Moreover, it is remarkable that the polynomials $(p_1x-1) \cdot (p_2x-1) \in \mathbb{Z}[x]$ possess a zero modulo $m$ for every $m \in \mathbb{N}$ but not a zero in $\mathbb{Z}$. The most famous polynomial fulfilling this property seems to be 
$$(x^2-13) \cdot (x^2-17) \cdot (x^2-13 \cdot 17) \in \mathbb{Z}[x]$$ from \cite[Ex. 1.0.4]{BS}, and indeed, the
above proof also works out (at least for completions $\mathbb{Z}_p$ instead of localizations $\mathbb{Z}_{(p)}$) using this polynomial. \\
\\
While the surjectivity of a polynomial map imposes severe restrictions on its Jacobian determinant, it is worth to note that singularities do not prevent polynomial maps from being injective in general. We shall devote a final proposition to this phenomenon. \\
\\
\textbf{Proposition 5.2:} \textit{Let $m \in \mathbb{N}$ and define $f = x^m \in \mathbb{Z}[x]$. Then $f \colon \mathbb{Z}_p \to \mathbb{Z}_p$ is injective if and only if $\gcd \big( m,2(p-1) \big) = 1$.} \\
\\
\textit{Proof:} We aim to study equations of the form $x^m = a$ which are solvable over $\mathbb{Z}_p$, i.e., under the assumption that there exists $b \in \mathbb{Z}_p$ so that $b^m = a$. Without loss of generality we may assume that $p \nmid a$ and thus that $p \nmid b$. \\
First, we can restrict to the case where $m = q$ is an odd prime number, as $x \mapsto x^m$ is injective if and only if $x \mapsto x^q$ is injective for every prime $q \mid m$, which clearly fails for $q = 2$. \\
Further, we note that the induced map $f \colon \mathbb{Z}/p \mathbb{Z} \to \mathbb{Z}/p \mathbb{Z}$ is bijective if and only if $q \nmid p-1$. This can be proved by examining the kernel of the homomorphism $x \in \big( \mathbb{Z}/p \mathbb{Z} \big)^\times \mapsto x^q \in \big( \mathbb{Z}/p \mathbb{Z} \big)^\times$ of groups. \\
That means, if $q \mid p-1$, then we find $c \in \mathbb{Z}$ so that $f(c) \equiv a \text{ (mod } p \text{)}$ but $c \not\equiv b \text{ (mod } p \text{)}$. Since $f'(c) = q \cdot c^{q-1} \not\equiv 0 \text{ (mod } p \text{)}$, there exists $\widehat{c} \in \mathbb{Z}_p, \widehat{c} \neq b,$ so that $f(\widehat{c}) = a$ by Hensel's lifting lemma and hence $f \colon \mathbb{Z}_p \to \mathbb{Z}_p$ is not injective. \\
On the other hand, if $q \nmid p-1$, then any $c \in \mathbb{Z}_p$ so that $f(c) = a$ suffices $b \equiv c \text{ (mod } p \text{)}$. If further $p \neq q$, then $f'(b), f'(c) \not\equiv 0 \text{ (mod } p \text{)}$ and hence Hensel's lifting lemma implies that $b = c$. If $p = q$, one passes from $f$ to $\tilde{f} = \frac{1}{p^2} \cdot (f(b + p \cdot x)-a)$ and draws the same conclusion. \\
Consequently, $f \colon \mathbb{Z}_p \to \mathbb{Z}_p$ is injective if and only if $2 \nmid m$ and no prime factor of $m$ divides $p-1$, thus if and only if $\gcd \big( m,2(p-1) \big) = 1$. $\hfill \Box$ \\
\\
Supposing that $m$ is odd, Proposition 5.2 implies that the map $x \in \mathbb{Z}_p \mapsto x^m \in \mathbb{Z}_p$ is injective if $p \equiv 2 \text{ (mod } m \text{)}$, but not injective if $p \equiv 1 \text{ (mod } m \text{)}$. By the Dirichlet prime number theorem, there exist infinitely many primes of both types. 
Furthermore, this also yields a number theoretic proof of the well-known fact that the map $x \in \mathbb{Z} \mapsto x^m \in \mathbb{Z}$ is injective if and only if $m$ is odd.

\subsection*{Acknowledgements}
The author is grateful to Giancarlo Castellano, Herwig Hauser and Arno van den Essen for their supportive and fruitful remarks, and indebted to an anonymous referee for providing a counterexample to Proposition 2.1 et al. if the assumption of finite quotients is dropped.

\textit{Postal address: Lukas Prader, Kirchberg 41A, 53179 Bonn, Germany} \\
\\
\textit{E-mail address: lukas.prader@uni-bonn.de}

\end{document}